\documentclass[11pt]{amsart}
\usepackage{amsmath}

\input epsf.def
\input epsf

\numberwithin{equation}{section}

\theoremstyle{plain}

\newtheorem{Main}{Main Theorem}

\newtheorem{Prop}{Proposition}[section]
\newtheorem{Thm}[Prop]{Theorem}		
\newtheorem{Cor}[Prop]{Corollary}
\newtheorem{Lem}[Prop]{Lemma}
\newtheorem{Def}[Prop]{Definition}

\theoremstyle{remark}

\newcommand{\CBbb}{\mathbb C}
\newcommand{\RBbb}{\mathbb R}
\newcommand{\ZBbb}{\mathbb Z}
\newcommand{\PBbb}{\mathbb P}
\newcommand{\HBbb}{\mathbb H}

\newcommand{\End}{\mathop{\rm End}\nolimits}
\newcommand{\dist}{\mathop{\rm dist}\nolimits}

\newcommand{\tr}{\mathop{\rm Tr}\nolimits}

\newcommand{\lra}{\longrightarrow}

\newcommand{\Hopf}{\mathop{\rm Hopf}\nolimits}

\newcommand{\G}{\mathfrak G}

\newcommand{\SL}{SL(2,\CBbb)}
\newcommand{\M}{\mathcal{M}_{\text{\tiny Higgs}}}
\newcommand{\PL}{\mathcal{P}\mathcal{L}(\Gamma)}
\newcommand{\R}{\mathcal{R}(\Gamma)}
\newcommand{\D}{\mathcal{D}(\Gamma)}
\newcommand{\X}{\mathcal{X}(\Gamma)}
\newcommand{\PSL}{\mathcal{P}\mathcal{S}\mathcal{L}(\Gamma)}
\newcommand{\PMF}{\mathcal{P}\mathcal{M}\mathcal{F}(\Gamma)}

\newcommand{\simrightarrow}{\buildrel \sim\over\longrightarrow }

\begin{document}


\title[The Morgan-Shalen Compactification]
	{On the Morgan-Shalen Compactification of the $\SL$ Character Varieties of Surface
Groups}

\author[Daskalopoulos]{G.  Daskalopoulos}

\address{Department of Mathematics \\
		Brown University \\
		Providence,  RI  02912}

\thanks{G.D. supported in part by NSF grant DMS-9504297}

\email{daskal@gauss.math.brown.edu}
\author[Dostoglou]{S. Dostoglou}

\address{Mathematical Sciences  Building \\ University of Missouri
\\Columbia MO  65211}

\email{stamatis@math.missouri.edu}

\thanks{S.D. supported in part by the Research Board of the University of
Missouri and the Arts \& Science Travel Fund of the University of Missouri,
Columbia}

\author[Wentworth]{R.  Wentworth}

\address{Department of Mathematics \\
   University of California \\
   Irvine,  CA  92697}

\email{raw@math.uci.edu}

\thanks{R.W. supported in part by NSF grant DMS-9503635 and a Sloan Fellowship}
\date{September 14, 1998}


\begin{abstract} A gauge theoretic description of the Morgan-Shalen compactification of the $\SL$
character variety of the fundamental group of a hyperbolic surface is given in terms of a natural
compactification of the moduli space of Higgs bundles via the Hitchin map.
\end{abstract}

\maketitle

\baselineskip=16pt

\section{Introduction}

Let $\Sigma$ be a closed, compact, oriented surface of genus $g\geq 2$ and  fundamental group
$\Gamma$.  Let
$\X$ denote the $\SL$ character variety of $\Gamma$, and $\D\subset \X$ the
closed subset consisting of conjugacy classes of discrete, faithful representations.  Then
$\X$ is an affine algebraic variety admitting a compactification $\overline{\X}$,
due to Morgan and Shalen \cite{MS1}, whose boundary points $\partial
\X=\overline{\X}\setminus \X$ correspond to elements of $\PL$, the space of
projective classes of length functions on $\Gamma$ with the weak topology.

Choose a metric $\sigma$ on $\Sigma$, and let $\M(\sigma)$ denote the moduli space of semistable
rank two Higgs pairs on $\Sigma(\sigma)$ with trivial determinant, as constructed by Hitchin
\cite{H}.  Then $\M(\sigma)$ is an algebraic variety, depending on the complex structure defined
by 
$\sigma$ (cf.\
\cite{Si}). By the theorem of Donaldson \cite{D}, $\M(\sigma)$ is homeomorphic to
$\X$, though not complex analytically so.  Let us denote this map $h:
\X\rightarrow\M$ (we henceforth assume the choice of basepoint $\sigma$).

We define a compactification of $\M$ as follows:  let $QD$ (more precisely, $QD(\sigma)$) denote
the finite dimensional complex vector space of holomorphic quadratic differentials on
$\Sigma$.  Then there is a surjective,
holomorphic map
$\M\to QD$ taking the Higgs field $\Phi$ to $\varphi=\det\Phi$.  We compose this with the map 
$$
\varphi\lra\frac{4\varphi}{1+4\Vert\varphi\Vert}
$$ where $\Vert\varphi\Vert=\int_\Sigma|\varphi|$, and obtain
$$
\widetilde\det : \M\lra BQD =\left\{\varphi\in QD : \Vert \varphi\Vert < 1\right\}
$$ Let $SQD=\left\{\varphi\in QD : \Vert \varphi\Vert = 1\right\}$ be the space of normalized
holomorphic quadratic differentials.  We then define
$\overline\M=\M\cup SQD$ with the topology given via the map $\widetilde\det$. The aim of this
paper is to compare the two compactifications $\overline{\X}$ and  $\overline\M$.

The points of $\PL$ may be regarded as arising from the translation lengths of minimal,
non-trivial
$\Gamma$ actions on
$\RBbb$-trees. Modulo isometries and scalings,   this correspondence is one-to-one, at least in
the non-abelian case (cf.\ \cite{CM} and Section \ref{S:Definitions} below). The boundary $\partial
\D$ consists of \emph{small actions}, i.e. those for which the arc stabilizer subgroups are all
cyclic.   With our choice of conformal structure $\sigma$ we can define  a continuous,
surjective map
\begin{equation}   \label{E:Hopf} H : \PL \lra SQD
\end{equation} 
When the length function $[\ell]$ is realized by the translation length function of  a  tree
\emph{dual} to the lift of a normalized holomorphic quadratic differential
$\varphi$, then $H([\ell])=\varphi$;  the full map is a continuous extension of this (see Theorem
\ref{T:Hcontinuity}) with the fibers of $H$ corresponding more generally to \emph{foldings} of dual
trees.

Let $\PMF$ denote the space of projective classes of measured foliations on $\Sigma$, modulo
isotopy and Whitehead equivalence.  By the theorem of Hubbard-Masur \cite{HM} we also have a
homeomorphism $HM : \PMF\simrightarrow SQD$.  It is not clear how to lift $H$ to
factor through $\PMF$ in a manner independent of $\sigma$.  However, it follows
essentially by Skora's theorem \cite{Sk} that if $H$ is restricted to
$\PSL$, the small actions, then it factors through $HM$ by a homeomorphism $
\PSL\simrightarrow\PMF$.

With this understood, we define a (set theoretic) map 
\begin{equation}  \label{E:hbar}
\bar h : \overline{\X}\lra
\overline\M
\end{equation}
 by extending the map $h$ to $H$ on the boundary.  We shall prove the following:
 
\begin{Main} The map $\bar h$ is continuous and surjective.  Restricted to the compactification
of the discrete, faithful representations $\overline{\D}$, it is a homeomorphism onto its
image.
\end{Main}

Note that the second statement follows from the first, since  $\partial \D$ consists of
small actions, and therefore the restricted map is injective by the above mentioned theorem of
Skora.  It would be interesting to determine the fibers of $\bar h$ in general; this
question will be taken up elsewhere.   We also remark that the
$SL(2,\RBbb)$ version of the above theorem leads to a harmonic maps  description of the Thurston
compactification of Teichm\"uller space and was first proved by Wolf \cite{W1}. 
 Generalizing this result to $\SL$ was one of the motivations for this paper.

This paper is organized as follows:  in Section \ref{S:Definitions} we review the Morgan-Shalen
compactification, the definition of the Higgs moduli space, and the notion of a harmonic map to an
$\RBbb$-tree.  In Section \ref{S:Trees}, we define the boundary map $H$.  The key point is that
the non-uniqueness in the correspondence between abelian length functions and $\RBbb$-trees
alluded to above nevertheless leads, via harmonic maps, to a well-defined geometric object on
$\Sigma$, in this case a quadratic
differential.  The most important result here is Theorem \ref{T:abelian}.  Along the way, we
give a criterion, Theorem \ref{T:uniqueness}, for uniqueness of harmonic maps to trees, using
the arguments in \cite{W3}.  The Main Theorem is then proven in Section
\ref{S:MethodI} as a consequence of our previous work \cite{DDW}. In the last section,
 a somewhat more concrete analysis of the behavior of high energy harmonic maps is
outlined, illustrating previous ideas.

\section{Definitions}       \label{S:Definitions}

Let $\Gamma$ be a hyperbolic surface group as in the introduction.  We denote by $\R$ the set
of representations of $\Gamma$ into $\SL$, and by $\X$ the set of characters of
representations.  Recall that a representation $\rho:\Gamma\to\SL$ defines a character
$\chi_\rho:\Gamma\to \CBbb$ by $\chi_\rho(g)=\tr\rho(g)$.  Two representations $\rho$ and $\rho'$
are \emph{equivalent} if $\chi_\rho=\chi_{\rho'}$.  It is easily seen (cf.\ \cite{CS}) that
equivalent irreducible representations are conjugate.  If $\rho$ is a reducible representation,
then we can write
$$
\rho(g)=\left(\begin{matrix} \lambda_\rho(g) & a(g) \\ 0 & \lambda_\rho(g)^{-1} \end{matrix}\right)
$$
for a representation $\lambda_\rho:\Gamma\to\CBbb^\ast$.  The character $\chi_\rho$ determines
$\lambda_\rho$ up to the inversion coming from the action of the Weyl group, 
and is in turn completely determined by it.  It is shown in \cite{CS} that the set of characters
$\X$ has the structure of  an affine algebraic variety.   

In \cite{MS1}, a (non-algebraic)
compactification $\overline\X$ of $\X$ is defined as follows:  let $C$ be the set of conjugacy
classes of $\Gamma$, and let $\PBbb(C)=\PBbb(\RBbb^C)$ be the (real) projective space of non-zero,
positive functions on $C$.  Define the map $\vartheta:\X\to\PBbb(C)$  by 
$$
\vartheta(\rho)=\left\{ \log\left(|\chi_\rho(\gamma)|+2\right) : \gamma\in C\right\}
$$
and let $\X^+$ denote the one point compactification of $\X$ with the inclusion map
$\imath: \X\to\X^+$.  Finally,  $\overline \X$ is defined to be the closure of the embedded
image of $\X$ in $\X^+\times\PBbb(C)$ by the map $\imath\times\vartheta$.  It is proved in
\cite{MS1} that $\overline \X$ is compact and that the boundary points consist of \emph{projective
length functions} on
$\Gamma$ (see the definition below).  Note that in its definition, $\vartheta(\rho)$ could be
replaced by the function $\left\{\ell_\rho(\gamma)\right\}_{\gamma\in C}$, where $\ell_\rho$
denotes the translation length for the action of $\rho(\gamma)$ on $\HBbb^3$:
$$
\ell_\rho(\gamma)=\inf\left\{\dist_{\HBbb^3}(x,\rho(\gamma)x) : x\in \HBbb^3\right\}
$$
(see \cite{Cp}).

Recall that an $\RBbb$-\emph{tree} is a metric space $(T,d_T)$ such that any two points
$x,y\in T$ are connected by a  \emph{segment} $[x,y]$, i.e.\ a rectifiable arc isometric to a
compact (possibly degenerate) interval in $\RBbb$ whose length realizes $d_T(x,y)$, and that
$[x,y]$ is the unique embedded path from
$x$ to $y$.  We say that
$x\in T$ is an
\emph{edge point} (resp.\ \emph{vertex}) if $T\setminus\{x\}$ has  two
(resp.\ more than two) components.    A $\Gamma$-\emph{tree} is an $\RBbb$-tree with an action of
$\Gamma$ by isometries, and it is called \emph{minimal} if there is no
proper $\Gamma$-invariant subtree.  We say that $\Gamma$ \emph{fixes an end} of $T$ (or more
simply, that $T$ has a fixed end) if there is a ray $R\subset T$ such that for every
$\gamma\in\Gamma$, $\gamma(R)\cap R$ is a subray.  When the action is understood, we shall often
refer to ``trees" instead of ``$\Gamma$-trees".

Given an $\RBbb$-tree $(T,d_T)$, the associated length function
$\ell_T:\Gamma\to\RBbb^+$ is defined by $\ell_T(\gamma)=\inf_{x\in T}d_T(x,\gamma
x)$.  If $\ell_T\not\equiv 0$, which is equivalent to $\Gamma$ having no fixed point
in $T$ (cf.\ \cite[Prop.\ II.2.15]{MS1}, then the class of $\ell_T$ in
$\PBbb(C)$ is called a projective length function.  We denote by $\PL$ the
set of all projective length functions on $\Gamma$-trees.
A length function is called \emph{abelian} if it is given by $|\mu(\gamma)|$ for
some homomorphism $\mu:\Gamma\to \RBbb$.  We shall use the following result:

\begin{Thm}[cf.\ \cite{CM}, Cor.\ 2.3 and Thm.\ 3.7]    \label{T:cullermorgan}
Let $T$ be a minimal $\Gamma$-tree  with  non-trivial length function $\ell_T$. 
Then
$\ell_T$ is non-abelian if and only if $\Gamma$ acts without fixed ends.  Moreover, if $T'$ is any
other minimal $\Gamma$-tree with the same non-abelian length function, then there is a unique
equivariant isometry $T\simeq T'$.
\end{Thm}

\noindent It is a  fact that  abelian
length functions, in general, no longer determine a unique minimal $\Gamma$-tree up to isometry
(e.g.\ see
\cite[Example 3.9]{CM}), and this presents one of the main difficulties dealt with in this paper.

We now give a quick review of the theory of Higgs bundles on Riemann surfaces and
their relationship to representation varieties.  Let $\Sigma$, $\Gamma$ be as in the
introduction.  A \emph{Higgs pair} is a pair $(A,\Phi)$, where $A$ is an $SU(2)$
connection on a rank $2$ smooth vector bundle $E$ over $\Sigma$, and $\Phi\in
\Omega^{1,0}\left(\Sigma,\End_0(E)\right)$, where $\End_0(E)$  denotes the bundle of
traceless endomorphisms of $E$.  The Hitchin equations are:
\begin{equation}   \label{E:hitchinequations}
\begin{aligned}
F_A+[\Phi,\Phi^\ast]&= 0 \\
D_A^{\prime\prime}\Phi &= 0
\end{aligned}
\end{equation}
The group $\G$ of (real) gauge transformations acts on the space of Higgs pairs and
preserves the set of solutions to (\ref{E:hitchinequations}).   We denote by $\M$ the set
of gauge equivalence classes of these solutions.  Then $\M$ is a complex analytic
variety of dimension $6g-6$ (the holomorphic structure depending upon the choice $\sigma$ on
$\Sigma$), which admits a holomorphic map (cf.\ \cite{H})
\begin{equation}  \label{E:hitchinmap}
\det:\M \lra QD=H^0(\Sigma, K_\Sigma^{\otimes 2}) \ :\
(A,\Phi) \mapsto \det\Phi=-\tr \Phi^2
\end{equation}
By associating to $[(A,\Phi)]\in \M$ the character of the flat $\SL$ connection $
A+\Phi+\Phi^\ast$, one obtains a homeomorphism (cf.\ \cite{D,C}) $h:\M\to\X$. 
Implicit in the definition of $h$ is a $\Gamma$-equivariant harmonic map $u$ from the
universal cover $\HBbb^2$ of $\Sigma$ to $\HBbb^3$.  It is easily verified
that the \emph{Hopf differential} of $u$, $\Hopf(u)=\tilde\varphi=\langle u_z, u_z\rangle dz^2$,
descends to a holomorphic quadratic differential $\varphi$ on $\Sigma$ equal to $\det \Phi$ (up to
a universal non-zero constant).

Having introduced harmonic maps, we now give an alternative way to view the Morgan-Shalen
compactification.  First, it follows by an easy application of the Bochner-Weitzenb\"ock formula
that a sequence of representations
$\rho_i$ diverges to the boundary  only if the energies $E(u_{\rho_i})$ of the associated
equivariant harmonic maps
$u_{\rho_i}$ are unbounded.  Furthermore, given such a sequence it was shown in 
\cite{DDW} that if the $\rho_i$ converge to a boundary point in the sense of
Morgan-Shalen, then the harmonic maps $u_{\rho_i}$ converge (perhaps after passing
to a subsequence) in the sense of Korevaar-Schoen to a $\Gamma$-equivariant harmonic
map $u:\HBbb^2\to(T,d_T)$, where $(T,d_T)$ is a minimal $\Gamma$-tree having the same
projective length function as the Morgan-Shalen limit of the $\rho_i$.  As pointed
out before, the tree is not necessarily uniquely defined, and even in the case where
the tree is unique, uniqueness of the harmonic map is problematic.  

Recall that a harmonic map to a tree means, by
definition, an energy minimizer for the energy functional defined in \cite{KS1}. 
Given such a map, its Hopf differential $\tilde\varphi$ can be defined almost everywhere,
and by
\cite[Lemma 1.1]{S1}, which can be adapted to the singular case, one can show that the
harmonicity of $u$ implies that $\tilde\varphi$ is a holomorphic quadratic differential.  The
equivariance of $u$ implies that $\tilde\varphi$ is the lift of a differential on $\Sigma$.  Note
also that if $u:\HBbb^2\to T$ is harmonic, then $\Hopf(u)\equiv 0$ if and only if $u$ is
constant.  In the equivariant case, this in turn is equivalent to $\ell_T\equiv 0$ (cf.\
\cite{DDW}).  For the rest of the paper, we shall tacitly assume $\ell_T\not\equiv 0$.

A particular example is the following: consider a non-zero holomorphic quadratic differential
$\varphi$, and denote by $\tilde\varphi$ its lift to $\HBbb^2$.  Let
$T_{\tilde\varphi}$ denote the vertical leaf space of the (singular) foliation of
$\tilde\varphi$, and let $\pi:\HBbb^2\to T_{\tilde\varphi}$ denote the natural projection. 
According to
\cite{MS2} (and using the correspondence between measured foliations and geodesic
laminations),
$T_{\tilde\varphi}$  is an
$\RBbb$-tree with an action of
$\Gamma$, and the projection $\pi$ is a $\Gamma$-equivariant continuous map.  We note two
important facts: (1) the vertices of $T_{\tilde\varphi}$  are precisely the image by $\pi$ of the
zeros of $\tilde\varphi$; and (2) since the action of $\Gamma$ on $T_{\tilde\varphi}$ is small,
$T_{\tilde\varphi}$  has no fixed ends (cf.\ \cite{MO}).

\begin{Prop}  \label{P:projection}
The map $\pi:\HBbb^2\to T_{\tilde\varphi}$ is harmonic with Hopf differential
$\tilde\varphi$.  
\end{Prop}

\begin{proof}  Since $T_{\tilde\varphi}$ has no fixed ends the
 existence of a harmonic map follows from \cite[Cor.\ 2.3.2]{KS2}.   The fact that $\pi$ is itself
an energy minimizer seems to be well-known (cf.\ the introduction to
\cite{GS}, and \cite{W2}): although the definition of harmonic map in the latter reference is
\`a priori different from the notion of an energy minimizer, a proof  follows
easily, for example, from the result in
\cite{W2}. Indeed, consider a sequence of
$\Gamma$-equivariant harmonic diffeomorphisms $u_i:\HBbb^2\to\HBbb^2$ with Hopf differentials
$t_i\tilde\varphi$, $t_i\to\infty$.  Let $d_i$ denote the pull-back distance functions on
$\HBbb^2$ by the $u_i$, and let $d_\infty$ denote the pseudo-metric obtained by pulling back the
metric on $T_{\tilde\varphi}$ by the projection $\pi$.  Extend all of these to pseudo-metrics, also
denoted
$d_{i}$ and $d_\infty$,  on the space $\HBbb^2_\infty$ constructed in  \cite{KS2}.  Then the
natural projection
$\HBbb^2\to\HBbb^2_\infty/d_\infty\simeq T_{\tilde\varphi}$ coincides with the map $\pi$.  On
the other hand, by \cite[Section 4.2]{W2}, $d_{i}\to d_\infty$ pointwise, locally
uniformly.  Therefore, by \cite[Thm.\ 3.9]{KS2},
$\pi$ is an energy minimizer.
\end{proof}

Next, we consider $\Gamma$-trees which are not necessarily of the form $T_{\tilde\varphi}$.  We
need the following:

\begin{Def} \label{D:morphism}
A \emph{morphism} of $\RBbb$-trees is a map $f:T\to T'$ such that every non-degenerate
segment $[x,y]$ has a non-degenerate subsegment $[x,w]$ such that  $f$
restricted to ${[x,w]}$ is a isometry onto its image.  The morphism $f$ is said to
\emph{fold} at a point $x\in T$ if there are non-degenerate segments $[x,y_1]$ and $[x,y_2]$
with $[x,y_1]\cap[x,y_2]=\{x\}$ such that $f$ maps each segment $[x,y_i]$ isometrically
onto a common segment in $T'$.  
\end{Def}

\noindent  It is a fact that a morphism $f:T\to T'$ is an isometric embedding
unless it folds at some point (cf.\ \cite[Lemma I.1.1]{MO}).  We also note that in general,
foldings $T\to T'$ may take vertices to edge points.  Conversely, vertices in $T'$
need not lie in the image of the vertex set of $T$.

\begin{Prop}[cf.\ \cite{FW}]  \label{P:factorization}
Let $T$ be an $\RBbb$-tree with $\Gamma$ action and let  $u:\HBbb^2\to T$ be an equivariant
harmonic map with Hopf differential $\tilde\varphi$.  Then $u$ factors as $u=f\circ \pi$,
where $\pi:\HBbb^2\to T_{\tilde\varphi}$ is as in Proposition \ref{P:projection} and
$f:T_{\tilde\varphi}\to T$ is an equivariant morphism.
\end{Prop}

\begin{proof}
Consider $f=u\circ \pi^{-1}:T_{\tilde\varphi}\to T$.  We first show that $f$ is
well-defined:  indeed, assume $z_1, z_2\in \pi^{-1}(w)$. Then $z_1$
and $z_2$ may be connected by a vertical leaf $e$ of the foliation of
$\tilde\varphi$.  Now by the argument in \cite[p.\ 117]{W3}, $u$ must collapse $e$ to
a point, so $u(z_1)=u(z_2)$.  In order to show that $f$ is a morphism, consider a segment $[x,z]\in
T_{\tilde\varphi}$.  We may lift $x$ to a point $\tilde x$ away from the zeros of
$\tilde\varphi$.  Moreover, we may choose a small horizontal arc $\tilde e$ from $\tilde x$
to some $\tilde y$ projecting to $[x,y]\subset[x,z]$, still bounded away from the zeros. 
The analysis in \cite{W3} again shows that this must map by $u$ isometrically onto
a segment in $T$.
\end{proof}

\noindent {\bf Remark.} It is easily shown (cf.\ \cite{DDW}) that images of equivariant harmonic
maps to trees are always minimal subtrees;  hence, throughout this paper we shall assume, without
loss of generality, that our trees are minimal.  Thus, for example, the factorization
$f:T_{\tilde\varphi}\to T$ above either folds at some point or is an equivariant isometry.

\section{The map $H$}          \label{S:Trees}

The basic fact is  that the Hopf differential for a harmonic map to a given tree is uniquely
determined:

\begin{Prop}   \label{P:Hopf} Let $T$ be a minimal $\RBbb$-tree with a non-trivial $\Gamma$
action.  If $u,v$ are equivariant harmonic maps $\HBbb^2\to T$, then $\Hopf(u)=\Hopf(v)$.
\end{Prop}

\begin{proof} This is proven in \cite{KS1}, where in fact the full pull-back ``metric tensor"
is considered.  In our situation, the result can also be seen as a direct consequence of the leaf
structure of the Hopf differential.   First, by \cite[p.\ 633]{KS1} the function
$z\mapsto d_T^2(u(z),v(z))$ is subharmonic; hence by the equivariance it must be equal to a
constant
$c$. We assume $c\neq 0$, since otherwise there is nothing to prove.  Set
$\tilde\varphi=\Hopf(u)$, $\tilde\psi=\Hopf(v)$. Suppose that
$p\in
\HBbb^2$ is a zero of
$\tilde\varphi$, and let $\Delta$ be a small neighborhood of $p$ containing no other zeros of
$\tilde\varphi$, and no zeros of $\tilde\psi$, except perhaps $p$ itself.  Then by Proposition
\ref{P:factorization} it follows that $u$ is constant equal to $u(p)$ on  every arc
$e\subset\Delta$ of the vertical foliation of
$\tilde\varphi$ with endpoint $p$.  On the other hand, $v(e)$ is a connected set satisfying
$d_T(u(p),v(z))=c$ for all $z\in e$.  Since spheres are discrete in trees, $v$ is constant equal
to $v(p)$ on $e$ as well.  Referring again to Proposition \ref{P:factorization}, this implies
that $e$ must be contained in a vertical leaf of $\tilde\psi$.  In this way, one sees that the
zeros of $\tilde\varphi$ and $\tilde\psi$ coincide with multiplicity in $\HBbb^2$. Thus, the same
is true for $\varphi$ and $\psi$ on
$\Sigma$.  Since the quadratic differentials are both normalized, they must be equal.
\end{proof}

We shall also need the following restriction on the kinds of foldings that arise from harmonic
maps:

\begin{Lem}  \label{L:noadjacentedges} Let $T_{\tilde\varphi}\to T$ arise from a harmonic map as
in Proposition \ref{P:factorization}.  Then folding occurs only at vertices, i.e.\ the images of
zeros of $\tilde\varphi$, and at the zeros of
$\tilde\varphi$, adjacent edges may not be folded.  In particular, folding cannot occur at simple
zeros.
\end{Lem} 

\begin{proof} The argument is similar to that in \cite[p.\ 587]{W2}. Suppose $p\in \HBbb^2$ is a
zero at which a folding occurs, and choose a neighborhood $\Delta$ of $p$ contained in a
fundamental domain and containing no other zeros.  We can find distinct segments
$e, e'$  of the horizontal foliation of
$\tilde\varphi$ with a common endpoint $p$  which map to segments of $T_{\tilde\varphi}$.  We
may further assume that the folding
$T_{\tilde\varphi}\to T$ carries each of $e$ and $e'$ isometrically onto a segment $\bar e$ of
$T$. Suppose that $e$ and $e'$ are adjacent. Then there is a small disk
$\Delta'\subset\HBbb^2$ which, under the projection $\pi:\HBbb^2\to T_{\tilde\varphi}$, maps to
$\pi(e)\cup \pi(e')$ and whose center maps  to $\pi(p)$ (see Fig.\  1).  Then the harmonic map
$u:\HBbb^2\to T$ maps
$\Delta'$ onto the segment
$\bar e$ with the center mapping to an endpoint.  Let $q$ denote the other endpoint of $\bar e$. 
The function $z\mapsto (d_T(u(z),q))^2$ is subharmonic on
$\Delta'$ with an interior maximum.  It therefore must be constant, which contradicts
$\varphi\not\equiv 0$. For the last statement, recall that the horizontal foliation is trivalent at
a simple zero, so that any two edges are adjacent.
\end{proof}

\bigskip
\centerline{\epsfbox{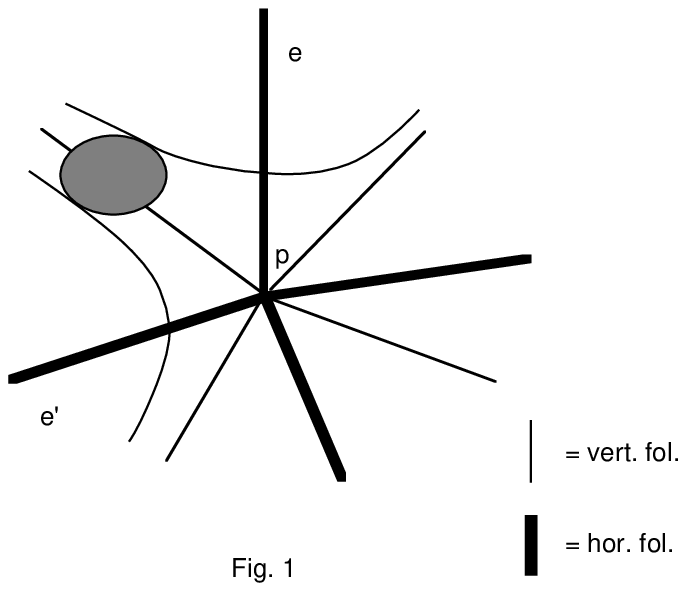}}
\bigskip

Though the following will not be important in this paper, it is interesting to note  that a
uniqueness result for equivariant harmonic maps to trees follows from these considerations, in
certain cases:

\begin{Thm}  \label{T:uniqueness}
 Let $u:\HBbb^2\to T$ be an equivariant harmonic map with 
$\tilde\varphi=\Hopf(u)$.  Suppose there is some vertex $x$ of $T_{\tilde\varphi}$ such that the
map
$f:T_{\tilde\varphi}\to T$ from Proposition \ref{P:factorization} does not fold at $x$.   Then
$u$ is the unique equivariant harmonic map to $T$.
\end{Thm}

\begin{proof}  Let $p$ be a zero of $\tilde\varphi$ projecting via $\pi$ to  $x$, and  let
$v$ be another equivariant harmonic map to $T$.  Choose a neighborhood $\Delta$ of $p$ as in the
proof of Proposition
\ref{P:Hopf}, and again suppose that the constant
$c=d_T(u(z),v(z))\neq 0$.  Recall that $x$ is a vertex of $T_{\tilde\varphi}$.  By the assumption
of no folding at
$x$, there must be a segment
$e$ of the vertical foliation of $\tilde\varphi$ in $\Delta$, with one endpoint being $p$, having
the following property:  for any $z\neq p$ in $e$ there is a neighborhood $\Delta'\subset\Delta$ of
$z$ such that
$u(\Delta')\cap [u(p),v(p)]=\{ u(p)\}$.  By Proposition
\ref{P:Hopf} and Lemma \ref{L:noadjacentedges}, we see that for
such $\Delta'$, $v(\Delta')\not\subset [u(p),v(p)]$.  Thus, there is a $q\in
\Delta$  such that $u(q)\not\in [u(p),v(p)]$ and $v(q)\not\in [u(p),v(p)]$. But then
$d_T(u(q),v(q))>d_T(u(p),v(p))=c$; contradiction.
\end{proof} 

\begin{Cor}  \label{C:uniqueness}
Let $\varphi\not\equiv 0$ be a holomorphic quadratic differential on $\Sigma$.
Then the map $\pi:\HBbb^2\to T_{\tilde\varphi}$ in Proposition \ref{P:projection} is the unique
equivariant harmonic map to $T_{\tilde\varphi}$.  If $u:\HBbb^2\to T$ is an equivariant harmonic
map and $\Hopf(u)$ has a zero of odd order, then $u$ is unique.
\end{Cor}

\begin{proof}  The first statement is clear from Theorem \ref{T:uniqueness}.  For the second
statement, notice that if $p$ is a zero of odd order we can still find a neighborhood $\Delta'$ as
in the proof of Theorem \ref{T:uniqueness}.
\end{proof}

Proposition \ref{P:Hopf} allows us to associate a unique $\varphi\in SQD$ to any non-abelian length
function:

\begin{Prop}  \label{P:nonabelian} 
 Let $[\ell]\in\PL$ be non-abelian.  Then there is a unique choice
 $\varphi\in SQD$ with the following property: if $T$ is any minimal $\RBbb$-tree with length
function $\ell$ in the class $[\ell]$, and
$u:\HBbb^2\to T$ is a $\Gamma$-equivariant harmonic map, then $\Hopf(u)=\varphi$.
\end{Prop}

\begin{proof} Let $\ell\in[\ell]$.  By Theorem \ref{T:cullermorgan}, there is a unique minimal tree
$T$, up to isometry, with length function $\ell$ and no fixed ends.  By Proposition \ref{P:Hopf},
any two harmonic maps
$u,v:\HBbb^2\to T$ have the same normalized  Hopf differential.  Furthermore, if $T'$ is
isometric to $T$ and $u'$ is a harmonic map to $T'$, then composing with the isometry, we see
that $u'$ has the same Hopf differential as any harmonic map to $T$.  If the length function
$\ell$ is scaled, then the normalized Hopf differential remains invariant.  Finally, since $T$
has no fixed ends, it follows from  \cite[Cor.\ 2.3.2]{KS2} that there exists an equivariant
harmonic map
$u:\HBbb^2 \to T$; so we set
$\varphi=\Hopf(u)$.
\end{proof}

We now turn our attention to the abelian length functions.  These no longer determine a unique
$\RBbb$-tree in general;  nevertheless, we shall see that there is still a uniquely defined
quadratic differential associated to them.  We begin with the following:

\begin{Prop}  \label{P:maptoreals} Let $\ell$ be an abelian length function, and let $\Gamma$ act
on $\RBbb$ with translation length function equal to $\ell$.  Then there is an equivariant
harmonic function $u:\HBbb^2\to\RBbb$, unique up to translations of $\RBbb$, with Hopf
differential $\tilde\varphi=(\tilde\omega)^2$, where
$\tilde\omega$ is the lift to $\HBbb^2$ of an abelian differential $\omega$ on $\Sigma$. 
 Moreover, $\ell$ is determined by the periods of $\text{Re}(\omega)$.
\end{Prop}

\begin{proof} The uniqueness statement is clear.  By harmonic theory there is a unique
holomorphic 1-form $\omega$ on
$\Sigma$ such that the real parts of its periods correspond to the homomorphism
$$
\mu : \pi_1(\Sigma)\lra H_1(\Sigma,\ZBbb)\lra \RBbb
$$ Choosing any base point $\ast$ of $\HBbb^2$, the desired equivariant harmonic function is the
real part of the holomorphic function $ f(z)=\int_\ast^z \tilde\omega$.  The Hopf
differential is $(f'(z))^2=(\tilde\omega)^2$.
\end{proof}

It is generally true that harmonic maps to trees with abelian length functions have Hopf
differentials with even order vanishing and that the length functions are recovered from the
periods of the associated abelian differential, as the next result demonstrates:

\begin{Thm}  \label{T:abelian} Let $u:\HBbb^2\to T$ be an equivariant harmonic map to a minimal
$\RBbb$-tree with non-trivial abelian length function $\ell$.  Then $\Hopf(u)=(\tilde\omega)^2$,
where $\tilde\omega$ is the lift to $\HBbb^2$ of an abelian differential $\omega$ on $\Sigma$. 
Moreover, $\ell$ is determined by the periods of $\text{Re}(\omega)$.
\end{Thm}

\begin{proof} We first prove that the Hopf differential $\tilde\varphi=\Hopf(u)$ must be a
square.  It suffices to prove that the zeros of $\tilde\varphi$ are all of even order.  Let $p$
be such a zero and choose a neighborhood $\Delta$ of $p$ as above.  Since $T$ has an abelian
length function, the action of
$\Gamma$ must fix an end $E$ of $T$.  Then applying the construction of Section 5 of \cite{DDW},
we find a continuous family of equivariant harmonic maps $u_{\varepsilon}$ obtained by ``pushing"
the image of $u$ a distance $\varepsilon$ in the direction of the fixed end. On the other hand, if
$\tilde\varphi$ had a zero of odd order, this would violate Corollary \ref{C:uniqueness}.

 We may
therefore express $\tilde\varphi=(\tilde\omega)^2$
 for some abelian differential
$\tilde\omega$ on $\HBbb^2$.  \`A priori, we can only conclude that $\tilde\omega$ descends to an
abelian differential $\hat\omega$ on an unramified double cover $\widehat \Sigma$ of $\Sigma$
determined by an index two subgroup $\widehat\Gamma\subset\Gamma$.  Let
$L$ be a complete non-critical leaf of the horizontal foliation of
$\tilde\varphi$.  Choose a point $x_0\in L$ and let $\bar x_0=u(x_0)$. We assume that we have
chosen $x_0$ so that $\bar x_0$ is an edge point. Then there is a unique ray
$\bar R$ with endpoint $\bar x_0$ leading out to the fixed end $E$.  Let $R$ denote the half-leaf
of $L$ starting at $x_0$ and such that a small neighborhood of $x_0$ in $R$ maps isometrically
onto a small subsegment of $\bar R$.

 We claim that $R$ itself  maps isometrically onto $\bar R$. For  suppose to the contrary that
there is a point $y\in R$ such that the portion $[x_0,y]$ of
$R$ from
$x$ to $y$ maps isometrically onto a subsegment of $\bar R$, but that this is not true for any
$y'\in R\setminus [x_0,y]$.  Clearly, the image of $y$ by $u$ must be a vertex of $T$.  
Recall the factorization $f:T_{\tilde\varphi}\to T$ from Proposition \ref{P:factorization}.  Since
$f$ is a surjective morphism of trees, the vertices of $T$ are either images by $f$ of vertices of
$T_{\tilde\varphi}$, and hence images by $u$ of zeros of $\tilde\varphi$, or they are vertices
created by a folding of $f$.  Thus, there are
two cases to consider: (1)  There is a point $q$ such that $y$ and $q$ lie on the same vertical
leaf and
$q$ is a zero of $\tilde\varphi$. Moreover, there is a critical  horizontal leaf $R'$ with one
endpoint equal to $q$, a small subsegment of which maps isometrically onto a subsegment of $\bar
R$ with endpoint $\bar q=u(q)$ (See Fig.\  2); and (2)  There is a point
$q$ such that $y$ and
$q$ lie on the same vertical leaf,
$q$ is connected by a horizontal leaf to a zero $p$ of $\tilde\varphi$, and the map
$f$ folds at $\pi(p)$, identifying the segment $[p,q]$ with a portion  $[p,q']$
of another horizontal leaf $R'$.  Moreover, $[p,q']$ maps isometrically onto a subsegment of the
unique ray from $\bar p =u(p)$ to the end $E$ (See Fig.\  3).  

\bigskip
\centerline{\epsfbox{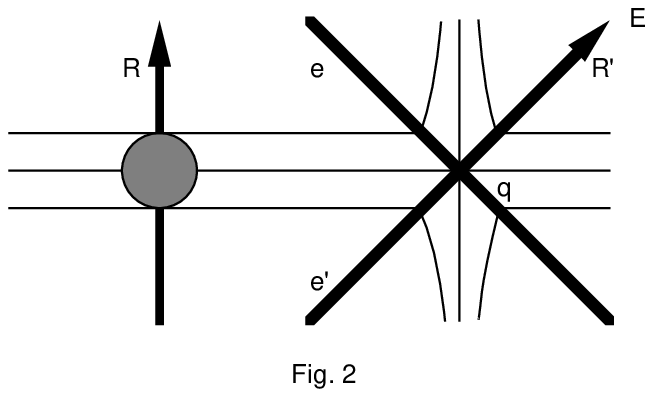}}
\bigskip

Consider Case (1):  As indicated in Fig.\  2, we can find a small neighborhood $\Delta$ of
$y$ and portions of horizontal leaves $e$ and $e'$ meeting at $q$ which map isometrically
onto segments of
$T$ intersecting the image $\bar R'=u(R')$ only in $\bar q$.  Now as above, by pushing the image
of $u$ in the direction of $E$, and possibly choosing $\Delta$ smaller, we can find a harmonic map
$u_\varepsilon$ which maps $\Delta$ onto a segment with endpoint $\bar q$, and maps $y$ to the
opposite endpoint; contradiction. The argument for Case (2) is similar:  We may find a disk
$\Delta$ centered at $y$ which maps to the union of segments $[\bar p,\bar q]$ and $[\bar r,\bar
q]$, with $y$ being mapped to $\bar q$.  Then pushing the map in the direction of $E$ as above
again leads to a contradiction (See Fig.\  3).

\bigskip
\centerline{\epsfbox{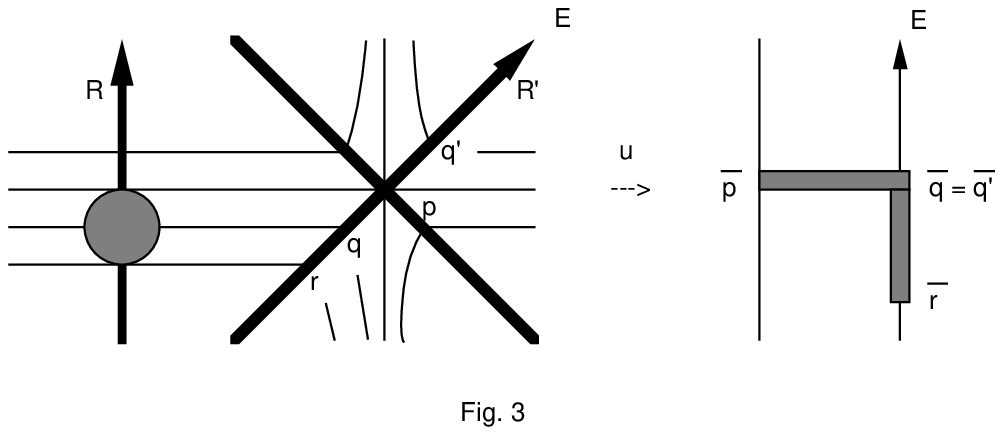}}
\bigskip

Next, we claim that for any $g\in \widehat\Gamma$, $\ell(g)$ is given by the period of
$\text{Re}(\hat\omega)$ around a curve representing the class $[g]$.  First, by
definition of a fixed end, the intersection  $\bar R\cap g(\bar
R)$   contains a subray of $\bar R$, and  for all $\bar x$ in this subray, 
$\ell(g)=d_T(\bar x,g(\bar x))$  (cf.\ \cite[Thm.\ 2.2]{CM}). For simplicity then, we assume
$g(\bar R)\subset\bar R$. Choose a lift of
$\bar x$ to
$x\in R$.  Then $u(g(x))=g(\bar x)\in \bar R$.  Suppose $g(x)$ is connected by a (possibly empty)
vertical leaf to a point $x'$ on $R$.  Then the curve $\tilde\gamma$ consisting of the portion
$[x,x']$ of $R$ from $x$ to $x'$ followed by the vertical leaf to $g(x)$ projects to a curve
$\gamma$ on $\Sigma$ representing $g$.  Moreover, since $R$ maps isometrically onto
$\bar R$, $\ell(g)$ is the length of $[x,x']$ with respect to the transverse measure
determined by $\tilde\varphi$.  Since $R$ contains no zeros of $\tilde\varphi$, the latter is
simply the absolute value of
$\int_{[x,x']}\text{Re}(\tilde\omega)$.  Futhermore, since the vertical direction
lies in the kernel of $\text{Re}(\hat\omega)$, we also have
$$
\ell(g)=\left|\int_{\tilde\gamma}\text{Re}(\tilde\omega)\right|=\left|
\int_\gamma\text{Re}(\hat\omega)\right|
$$ as desired.

Now consider the possibility that $g(x)\in g(R)$ is not connected to $R$ by a vertical leaf. 
Since $g(\bar x)\in\bar R$, it follows from Proposition \ref{P:factorization} and the fact that
$R$ maps onto $\bar R$ that there is an intervening folding of a subray of $g(R)$ onto $R$.  Let
$y\in R$ project to the vertex in
$T_{\tilde\varphi}$ at which this occurs.  The simplest case is where $y$ is connected by a
vertical leaf to a point  $w\in g(R)$, and the folding identifies the subray of $R$ starting at
$y$ isometrically with the subray of $g(R)$ starting at $w$.  The same analysis as above then
produces the closed curve $\gamma$.

A more complicated situation arises when there are intervening vertices (See Fig.\ 4(a)): 
For example, there may be zeros $p,q$  of $\tilde\varphi$, a point $w'\in g(R)$, and segments $e,
e', e''$ of the vertical, horizontal, and vertical foliations, respectively, with endpoints
$\{y,p\}$, $\{p,q\}$, and $\{q, w'\}$, respectively.  Moreover, the map $u$
folds
$e'$ onto a subsegment $f$ of $R$ with endpoints $y$ and $y'$, and then identifies the subray of
$R$ starting at $y'$ isometrically with the subray of $g(R)$ starting at $w'$.  In this way, we
see that a subsegment $f'$ of $g(R)$ with endpoints $w'$ and $w$ gets identified with $f$ and
$e'$; in particular, the transverse measures of these three segments are all equal (strictly
speaking,
$y'$ need not lie on
$R$ as we have chosen it, but this will not affect the argument).

Now consider the prongs at the zero $p$, for example.  These  project to distinct
segments in $T_{\tilde\varphi}$, which are then either projected to  segments in
$T$ intersecting
$\bar R$ only in $\bar y$, or alternatively there may be a folding identifying them with
subsegments of $\bar R$.  Let us label the prongs with a $+$ sign if there is a folding onto a
subsegment of $[\bar y, E)$, with a $-$ sign if there is a folding onto a subsegment of $[\bar
x,\bar y]$, and with a $0$ if no folding occurs, or if the edge is folded along some other segment
(See Fig.\ 4(b)).  Since
$p$ is connected by the vertical leaf $e$ to $R$, we label the adjacent horizontal segments with
$+$ and $-$ accordingly.  Working our way around $p$ in the clockwise direction, and repeatedly
using the ``pushing" argument from Section 5 of \cite{DDW}, we find that every second prong must be
labelled
$+$ while the intervening prongs may get either $-$ or $0$ (recall here Lemma
\ref{L:noadjacentedges}).  Therefore, there must be an odd number of prongs between $e'$ and the
one adjacent to $e$ which is identified in the leaf space with a portion of
$f$.  A similar argument applies to $q$, $e''$, and $f'$.

\bigskip
\centerline{\epsfbox{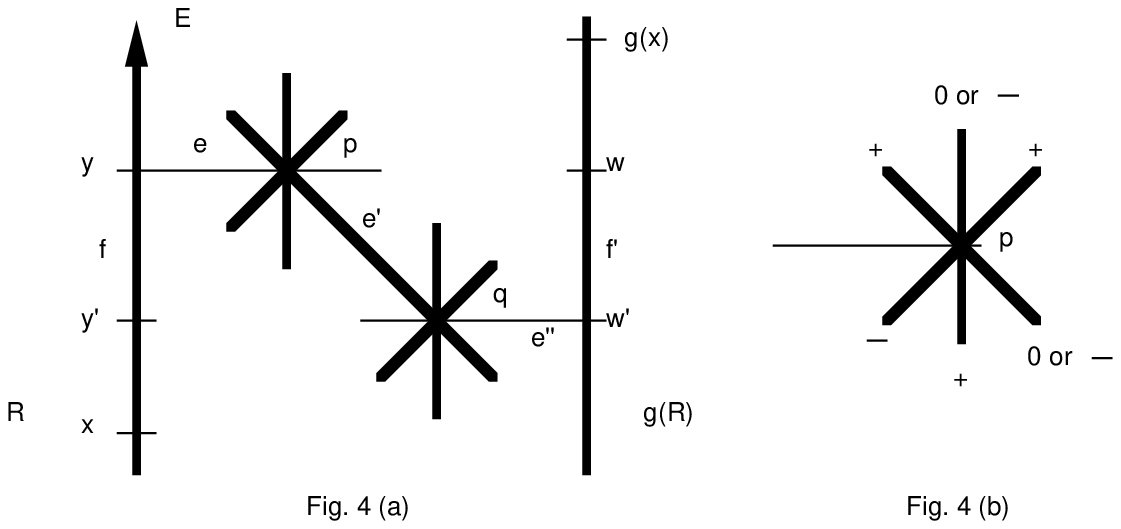}}
\bigskip

Let $\tilde\gamma'$ be the path from $y'$ to $w$ obtained by following $f, e, e', e''$, and then
$f'$.  Because of the odd sign to the folding of the prongs at $p$ and $q$, one may easily verify
that $\left|\int_{\tilde\gamma'}\text{Re}(\tilde\omega)\right|$ is the just the transverse
measure of the segment $f$. Indeed, suppose $\tilde\varphi$ has a zero of order $2n$ at some point
$p$, and choose a local conformal coordinate $z$ such that $\tilde\varphi(z)= z^{2n}dz^2$. Then
the foliation is determined by the leaves of $\xi=z^{n+1}/n+1$.  If $\zeta$ is  a primitive
$2n+2$ root of unity, then $z\mapsto \zeta^k z$ takes one radial prong to another, with $k-1$
prongs in-between (in the counter-clockwise direction).  The outward integrals of
$\text{Re}\sqrt{\tilde\varphi}$ along these prongs to a fixed radius differ by $(-1)^k$.  Our
analysis implies that $k-1$ is odd, so $k$ is even, and we have the correct cancellation.  
If we
extend
$\tilde\gamma'$ along the horizontal leaves
$R$ and
$g(R)$ to a path $\tilde\gamma$ from $x$ to $g(x)$, then
$\left|\int_{\tilde\gamma}\text{Re}(\tilde\omega)\right|=d_T(\bar x, g(\bar x))$ as required.  In
general, there will be additional intervening zeros, and the procedure above applies to each of
these with no further complication. 

Thus, $\ell$ restricted to $\widehat\Gamma$ is given by the periods of $\text{Re}(\hat\omega)$. 
Since the real parts of the periods of an abelian differential determine the differential
uniquely, $\hat\omega$ must agree with the pull-back to $\widehat\Sigma$ of the  form
in Proposition \ref{P:maptoreals}; in particular, it descends to $\Sigma$.
 This completes the proof of Theorem \ref{T:abelian}.
\end{proof}

We immediately have the following:

\begin{Cor}  \label{C:abelian}
Fix an abelian length function $\ell$.  Then for any tree $T$ with length function $\ell$ and
any equivariant harmonic map 
$v:\HBbb^2\to T$, we have $\Hopf(v)=\Hopf(u)$ where $u$ is the equivariant harmonic function
from Proposition \ref{P:maptoreals} corresponding to $\ell$.
\end{Cor}

We are now prepared to define the map (\ref{E:Hopf}).  Take a representative $\ell$ of
$[\ell]\in\PL$.  There are two cases:  if $\ell$ is non-abelian, use Proposition
\ref{P:nonabelian} to define $H([\ell])=\varphi$;  if $\ell$ is abelian, use Proposition 
\ref{P:maptoreals}.  The main result of this section is the following:

\begin{Thm}  \label{T:Hcontinuity} The map $H:\PL\to SQD$ defined above is continuous.
\end{Thm}

\begin{proof} Suppose $[\ell_i]\to [\ell]$, and assume, to the contrary, that there is a
subsequence, which we take to be the sequence itself, such that
$H([\ell_i])\to \varphi\neq H([\ell])$.  Choose representatives $\ell_i\to\ell$.  If there is a
subsequence $\{i'\}$ consisting entirely of abelian length functions, then $\ell$ itself must be
abelian, and from the construction of Proposition \ref{P:maptoreals},
$H(\ell_{i'})\to H(\ell)$; contradiction.  Thus, we may assume all the $\ell_i$'s are
non-abelian.  There exist $\RBbb$-trees $T_i$, unique up to isometry, and equivariant harmonic
maps
$u_i:\HBbb^2\to T_i$. We claim that the $u_i$ have uniform modulus of continuity (cf.\
\cite[Prop.\ 3.7]{KS2}).  Indeed, by \cite[Thm.\ 2.4]{GS}, it suffices to show that $E(u_i)$ is
uniformly bounded.  If $E(u_i)\to\infty$, then the same argument as in \cite[proof of Thm.\
3.1]{DDW} would give a contradiction.  It follows by \cite[Prop.\ 3.7]{KS2} that there is a
subsequence $\{i'\}$ (which we assume is the sequence itself) such that $u_i$ converges in the
pullback sense to an equivariant harmonic map $u:\HBbb^2\to T$, where $T$ is a minimal
$\RBbb$-tree with length function equal to $\ell$.  In addition, by \cite[Theorem 3.9]{KS2}, 
$\Hopf(u_i)\to\Hopf(u)$. If
$\ell$ is non-abelian, we have a contradiction by Proposition \ref{P:Hopf};  if $\ell$ is abelian,
we have a contradiction by Corollary
\ref{C:abelian}.
\end{proof}

\section{Proof of the Main Theorem}  \label{S:MethodI}

We show how the results of the previous section, combined with those in  \cite{KS2} and
\cite{DDW}, give a proof of the Main Theorem.  
We first reduce the proof of the continuity of $\bar h$ to the following:

\medskip\noindent 
{\bf Claim.}\  \emph{If $[\rho_i]\in \X$ is a sequence of representations
converging to $[\ell]\in\PL$ then $h([\rho_i])\to H([\ell])$.}

\medskip\noindent
For suppose the claim holds and $\bar h$ is not continuous.  Then we may find a sequence
$x_i\in\PL\cup \X$ such that $x_i\to x$ but $\bar h(x_i)\to y\neq \bar h(x)$.  If $x\in
\PL$ so that $\bar h(x)=H(x)$, the claim rules out the possibility that there is a subsequence of
$\{x_i\}$ in
$\X$.  In this case then, there must be a subsequence in $\PL$. But this contradicts the
continuity of $H$, Theorem \ref{T:Hcontinuity}.  Thus, $x$ must be in $\X$.  But then we
may assume that $\{x_i\}\subset \X$, so that $\bar h= h$ on $\{x_i\}$.  The continuity
of  the homeomorphism $h: \X\to\M$ then provides the contradiction.

It remains to prove the claim. 
Again suppose to the contrary that $[\rho_i]\to [\ell]$ but
$h([\rho_i])\to \varphi\neq H([\ell])$ for $\varphi\in SQD$.
 First, suppose that there is a subsequence $[\rho_{i'}]$
 with
reducible representative representations $\rho_{i'}:\Gamma\to\SL$.  Up to conjugation, 
which amounts to changing the choice of representative, we may assume each $\rho_{i'}$ fixes a
given vector
$0\neq  v\in \CBbb^2$, and that the action on the one dimensional line spanned by $v$ is
determined by a character $\chi_{i'}:\Gamma\to \CBbb^\ast$.  The associated translation length
functions
$\ell_{i'}$ are therefore all abelian, and so $[\ell]$ must be abelian.  We may assume there is a
representative $\ell$ such that $\ell_{i'}\to\ell$.
By Proposition \ref{P:maptoreals} there are 
harmonic functions
$$
u, u_{i'} : \HBbb^2\lra \RBbb\simeq\CBbb^\ast / U(1)\hookrightarrow\HBbb^3
$$
equivariant for the induced action of $\Gamma$ on $\CBbb^\ast$ by $\chi$ and $\chi_{i'}$,
respectively,  and these converge (after rescaling) to a harmonic function $u:\HBbb\to\RBbb$
equivariant with respect to an action on $\RBbb$ with translation length function $\ell$.  Since
the length functions converge, it follows from the construction in  Proposition \ref{P:maptoreals}
that
$\Hopf(u_{i'})\to\Hopf(u)$, so by the definition of $H$, $h([\rho_{i'}])\to H([\ell])$;
contradiction.

Second, suppose that there is a subsequence $[\rho_{i'}]$ of irreducibles.  Then by the main
result of \cite{DDW} we can find a further subsequence (which we take to be the sequence itself)
of $\rho_{i'}$-equivariant harmonic maps $u_{i'}:\HBbb^2\to\HBbb^3$ converging in the sense of
Korevaar-Schoen to a harmonic map $u:\HBbb^2\to T$, where $T$ is a minimal $\RBbb$-tree with an
action of $\Gamma$ by isometries and length function $\ell$ in the class $[\ell]$.  As above,
$\Hopf(u_{i'})\to\Hopf(u)$, so  by the definition of $H$, $h([\rho_{i'}])\to H([\ell])$;
contradiction.  Since we have accounted for both possible cases, this proves the claim.

\section{Convergence of Length Functions}  \label{S:MethodII}

In this last section we would like to briefly sketch an alternative argument for the convergence to
the boundary in the Main Theorem based on a direct analysis of length functions, more in the
spirit of
\cite{W1}.  The  generalization of estimates for equivariant harmonic maps with target $\HBbb^2$ 
to maps with target
$\HBbb^3$ has largely been carried out by Minsky \cite{M}.  We discuss this point of view,
however, since it reveals  how and why the folding of the dual tree $T_{\tilde\varphi}$
occurs.

The first step is to analyze the behavior of the induced metric for a harmonic map
$u:\HBbb^2\to\HBbb^3$ of high energy (at the points where $u$ is an immersion). As usual we will
denote by $\tilde\varphi$ the Hopf differential for the map $u$.  Because of equivariance,
$\tilde\varphi$ will be the lift of a holomorphic quadratic differential $\varphi$ on $\Sigma$; 
 recall the norm $\Vert \varphi\Vert$ from the introduction, and let $Z(\varphi)\subset \Sigma$
denote the zero set of $\varphi$.  We also  set
$\mu$ to be the Beltrami differential associated to the pull-back metric $u^\ast ds^2_{\HBbb^3}$. 

\begin{Lem}     \label{L:expdecay} Fix $\delta, T>0$.  Then there are constants
$B,\alpha>0$ such that for all $u, \mu$  and $\varphi$ as above, $\Vert\varphi\Vert\geq T$, and all
$p\in
\Sigma$ satisfying
$\dist_\sigma(p,Z(\varphi))\geq \delta$ we have
$$
\log\left(1/|\mu|\right)(p) < Be^{-\alpha \Vert\varphi\Vert}\ .
$$
\end{Lem}

\begin{proof}  This result is proven in \cite[Lemma 3.4]{M}.  One needs only a statement
concerning the uniformity of the constants appearing there.  However, by using the compactness of
$SQD$, one easily shows that for $\delta>0$ there is a constant $
c(\delta)>0$ such that  for all $\varphi\in SQD$ and all
$p\in \Sigma$ such that $\dist_\sigma(p,Z(\varphi))\geq \delta$ the disk
$U$ of radius $\tilde c(\delta)$ (with respect to the singular flat metric $|\varphi|$) around $p$
is embedded in
$\Sigma$ and contains no zeros of
$\varphi$.  Then the result cited above applies.
\end{proof}

This estimate is all that is needed to prove convergence in the case where there cannot be a
folding of the dual tree $T_{\tilde\varphi}$ such that the composition of projection to
$T_{\tilde\varphi}$ with the folding is harmonic.  From Lemma
\ref{L:noadjacentedges}, this will be guaranteed, for example,  if $\varphi$ has only simple zeros.
For simplicity, in this section we assume all representations are irreducible.

\begin{Thm}    \label{T:nofolding} 
Given an unbounded sequence $\rho_j$ of representations with
Morgan-Shalen limit 
$[\ell]$,
let $u_j:\HBbb^2\to\HBbb^3$ be the associated $\rho_j$-equivariant harmonic maps.
 Suppose that for $\tilde\varphi_j=\Hopf(u_j)$ we have
$\varphi_j/\Vert\varphi_j\Vert\to\varphi\in SQD$, where
$\varphi$ has only simple zeros.  Then $[\ell]=[\ell_T]$, where 
$T=T_{\tilde\varphi}$.  
\end{Thm}

\begin{proof}
We will prove the convergence of length functions in two steps.  First, we compare the length of
 closed curves $\gamma$ in the  free isotopy class $[\gamma]$ with respect to the induced metric
from
$u_j$ to the length with respect to the transverse measure.  Second, we will compare the length of
the image by $u_j$ of a lift
$\tilde\gamma$ to $\HBbb^2$ of $\gamma$  to the translation length in $\HBbb^3$
of the conjugacy class $[\gamma]$ represents.  The basic idea is that the image of $\tilde\gamma$ 
very nearly approximates a segment of the hyperbolic axis for $\rho_j([\gamma])$.

   For
$\varphi$ and
$[\gamma]$ as above, let
$\ell_\varphi([\gamma])$ denote the infimum over all representatives $\gamma$ of $[\gamma]$ of the
length of $\gamma$ with respect to the vertical measured foliation defined by $\varphi$.  If
$u:\HBbb^2\to\HBbb^3$ is a differentiable equivariant map, we define $\ell_u([\gamma])$ as
follows:  for each representative $\gamma$ of $[\gamma]$, where $[\gamma]$ corresponds to
 the conjugacy class of
$g\in\Gamma$, lift $\gamma$ to a curve $\tilde\gamma$ at a point $x\in \HBbb^2$, terminating at
$gx$.  We then take the infimum over all such $\tilde\gamma$ of the length of $u(\tilde\gamma)$. 
This is $\ell_u([\gamma])$, and by the equivariance of $u$ it is independent of the choice of $x$. 
Finally, recall that the translation length
$\ell_\rho([\gamma])$ for a representation $\rho:\Gamma\to\SL$ has been defined in Section
\ref{S:Definitions}.

Given $\varepsilon>0$, let $QD_\varepsilon\subset QD\setminus\{0\}$ denote the subset consisting of
holomorphic quadratic differentials $\varphi$ having only simple zeros, and such that the zeros
are pairwise at least a $\sigma$-distance $\varepsilon$ apart. Notice that for $t\, \neq 0$, $t\,
QD_\varepsilon=QD_\varepsilon$. The next result is a consequence of Lemma \ref{L:expdecay}:

\begin{Prop}     \label{P:pullbacklength} For all classes $[\gamma]$ and differentials
$\varphi\in QD_\varepsilon$ there exist constants $k$
 and
$\eta$ depending on $\Vert\varphi\Vert$, $[\gamma]$, and $\varepsilon$, so that
$$ k\, \ell_\varphi([\gamma])+\eta\geq \ell_u([\gamma])\geq \ell_\varphi([\gamma])
$$ where $k\to 1$ and $\eta\Vert\varphi\Vert^{-1/2}\to 0$ as $\Vert\varphi\Vert\to\infty$ in
$QD_\varepsilon$.
\end{Prop}

\begin{proof}[Sketch of proof]
We first need to choose an appropriate representative for the class of $[\gamma]$.  Such a choice
was explained in \cite{W1}.  Namely, for  $\delta>0$ and a given $\varphi$, we can find a
representative
$\gamma$ consisting of alternating vertical and horizontal segments and having the transverse
measure of the class
$[\gamma]$.  Moreover, because the zeros of  $\varphi$ are simple,
for sufficiently small $\delta$ we can also guarantee that $\gamma$ avoid a
$\delta$ neighborhood of the zeros.  
Now the proof follows as in \cite[Lemma 4.6]{W1}.  Note that along a \emph{harmonic maps ray},
i.e.\ a sequence $u_i$ such that $\Hopf(u_i)$ is of the form $t_i\varphi$ for a fixed $\varphi$
and an increasing unbounded sequence $t_i$, we no longer necessarily have monotonicity of the norm
of the Beltrami differentials $|\mu(t_i)|$.  The argument for the estimate still applies, however,
since
 the representatives $\gamma$ are uniformly supported away from the zeros.  There, we apply the
estimate Lemma \ref{L:expdecay}.  The details are omitted.
\end{proof}

Next, we compare $\ell_u$ with the translation length in $\HBbb^3$:

\begin{Prop}     \label{P:translationlength} 
Let $\rho:\Gamma\to\SL$ and $u:\HBbb^2\to\HBbb^3$ be the $\rho$-equivariant harmonic map with
$\tilde\varphi=\Hopf(u)$.  Suppose $\varphi\in QD_\varepsilon$.   For  all classes
$[\gamma]$ there exist constants $m$
 and
$\zeta$ depending on $\Vert\varphi\Vert$, $[\gamma]$, and $\varepsilon)$, so that
$$ m\, \ell_{\rho}([\gamma])+\zeta\geq \ell_{u}([\gamma])\geq \ell_{\rho}([\gamma])
$$ where $m\to 1$ and $\zeta\Vert\varphi\Vert^{-1/2}\to 0$ as
$\Vert \varphi\Vert \to\infty$ in $QD_\varepsilon$.
\end{Prop}

Combining Propositions \ref{P:pullbacklength} and \ref{P:translationlength} proves the theorem.
\end{proof}

\begin{proof}[Sketch of proof of Proposition \ref{P:translationlength}]
One observes that away from the zeros the images of the horizontal leaves of the foliation of
$\tilde\varphi$ are closely approximating (long) geodesics in $\HBbb^3$, while by Lemma
\ref{L:expdecay} the images of vertical leaves are collapsing.  More precisely, the following is
proven in
\cite[Thm.\ 3.5]{M}:

\begin{Lem}   \label{L:geodesic} 
Fix $\delta>0$, a representation $\rho: \Gamma\to \SL$, and let
$u:\HBbb^2\to\HBbb^3$ be the $\rho$-equivariant harmonic map with Hopf differential
$\tilde\varphi$. Let
$\tilde\beta$ be a segment of the horizontal foliation of $\tilde\varphi$ from $x$ to $y$ and
suppose that for all $\tilde p\in \tilde\beta$, $\dist_\sigma(p,Z(\varphi))\geq \delta$.  Then
there is an
$\varepsilon$, exponentially decaying in $\Vert\varphi\Vert$, such that
\begin{enumerate}
\item  $u(\tilde\beta)$ is uniformly within $\varepsilon$ of the geodesic in $\HBbb^3$ from $u(x)$
to
$u(y)$.
\item  The length of $u(\tilde\beta)$ is within $\varepsilon$ of $\dist_{\HBbb^3}(u(x),u(y))$.
\end{enumerate}
\end{Lem}

The following is the key result:

\begin{Lem}  \label{L:axis}
 Given $g\in \SL$,  let $\ell(g)$ denote the translation length for the action of $g$ on
$\HBbb^3$.  Suppose that
$s\subset\HBbb^3$ is a curve which is
$g$ invariant and satisfies the following property: for any two points $x,y\in s$, the segment of
$s$ from $x$ to $y$ is uniformly within a distance $1$ of the geodesic in $\HBbb^3$ joining $x$ and
$y$.  Then there is a universal constant $C$ such that 
$$\inf_{x\in s}\dist_{\HBbb^3}(x,g(x))\leq \ell(g)+C\ .$$
\end{Lem}

\begin{proof}  The intuition is clear: such an $s$ must be an ``approximate axis" for $g$.  The
proof proceeds as follows:  choose $x\in s$, and let $c$ denote the geodesic in $\HBbb^3$
from $x$ to $g(x)$.  By \cite[Lemma 2.4]{Cp} there exists a universal constant $D$ and  a
subgeodesic
$\tilde c$ of
$c$ with the property that $\left|\text{length}(\tilde c)-\ell(g)\right|\leq D$.  Let $a$
and $b$ be the endpoints of $\tilde c$ closest to $x$ and $g(x)$, respectively.  By the
construction of $\tilde c$ in the reference cited above, it follows that
$\dist_{\HBbb^3}(b,g(a))\leq D$; hence,
$\dist_{\HBbb^3}(b,g(b))\leq \ell(g)+2D$.  Now by the assumption on $c$, there is a point $y\in s$
close to $b$, so that $\dist_{\HBbb^3}(y,g(y))\leq \ell(g)+C$, where $C=2(D+1)$.
\end{proof}

Proceeding with the proof of Proposition \ref{P:translationlength}, choose the representative
$\gamma$ as discussed in Proposition
\ref{P:pullbacklength}.  We may then lift to $\tilde\gamma\subset \HBbb^2$ so that $\tilde\gamma$
is invariant under the action of
$g$.  Now $\tilde\gamma$ is written as a union of horizontal and vertical segments of the
foliation of $\tilde\varphi$.  Let $s=u(\tilde\gamma)$.  Then Lemmas \ref{L:expdecay} and
\ref{L:geodesic} imply that
$s$ satisfies the hypothesis of Lemma \ref{L:axis}.  Moreover, using Lemma \ref{L:geodesic} again,
along with some elementary hyperbolic geometry, one can show that $\inf_{x\in
s}\dist_{\HBbb^3}(x,g(x))$ is approximated by the length of a segment of $u(\tilde\gamma)$ from a
point $u(x)$ to $u(gx)$.   We leave the precise estimates to the reader.
\end{proof}

From Lemma \ref{L:noadjacentedges}, we see that foldings can only arise when the Hopf
differentials converge in $SQD$ to differentials with multiplicity at the  zeros. From the
point of view taken here, this corresponds to the fact that the representatives for closed curves
$\gamma$ chosen above may be forced to run into  zeros of the Hopf differential where the estimate
Lemma
\ref{L:expdecay} fails.  These may cause non-trivial angles to form in the image $u(\tilde\gamma)$
which, in the limit, may fold the dual tree.  

 Consider again the situation along a harmonic maps
ray with differential $\varphi$.  Given $[\gamma]$ corresponding to the conjugacy class of an
element $g\in\Gamma$, representatives
$\gamma$ still may be chosen as in the proof of Proposition \ref{P:pullbacklength} so that the
horizontal segments remain bounded away from the zeros.  However, it may happen that a vertical
segment passes through a zero of order
$2$ or greater.  For simplicity, assume this happens once.  Divide $\gamma$ into curves $\gamma_1$,
$\gamma_2$, and
$\gamma_v$, where $\gamma_v$ is the offending vertical segment, and lift to segments
$\tilde\gamma_1$,
$\tilde\gamma_2$, and
$\tilde\gamma_v$ in $\HBbb^2$.  Note that one endpoint of each of the $\tilde\gamma_i$'s
corresponds to either endpoint of $\tilde\gamma_v$, and the other endpoints of the
$\tilde\gamma_i$'s are related by
$g$.   By the Lipschitz estimate for harmonic maps to non-positively curved spaces we have a bound
on the distance in
$\HBbb^3$ between the endpoints of $u(\tilde\gamma_v)$  in terms of the length of $\gamma_v$ and
the energy
$E(u)^{1/2}$ (cf.\ \cite{S2}).  Thus, the rescaled length is small; in fact, since the length of
$\gamma_v$ is arbitrary, the distance converges to zero.  On the other hand, the previous argument
applies to the segments $u(\tilde\gamma_1)$ and $u(\tilde\gamma_2)$ which are connected by
$u(\tilde\gamma_v)$.  Adding the geodesic in $\HBbb^3$ joining the other endpoints of
$u(\tilde\gamma_1)$ and
$u(\tilde\gamma_2)$ forms an approximate geodesic quadrilateral which, in the rescaled limit,
converges either to an edge $\mid$ (no folding) or a possibly degenerate tripod $\dashv$
(folding).  In both cases, there is an edge which, by the same argument as in the
proof of Proposition
\ref{P:translationlength}, approximates the axis of $\rho_j(g)$ for large $j$.  At the same time,
the rescaled length of this segment is approximated by the translation length of the element $g$
acting on a folding of $T_{\tilde\varphi}$ at the zero.

An interesting question is whether this approach may be used to determine precisely the fibers of
the map $\bar h$ in the Main Theorem.  While the essential ideas have been outlined here, a
complete description is not yet available.  We will return to this issue in a future work.

\end{document}